\documentclass[10pt]{amsart}
\usepackage{amssymb}
\usepackage{dsfont}
\usepackage{amscd}
\usepackage[mathscr]{euscript} 
\usepackage[all]{xy}
\usepackage{xcolor}
\usepackage{url}
\usepackage{stmaryrd}
\usepackage{comment}
\usepackage[retainorgcmds]{IEEEtrantools}
\usepackage{hyperref}
\usepackage{graphicx}
\usepackage{mathtools}
\usepackage{centernot}
\usepackage{marvosym}
\usepackage{tikz}
\usepackage{soul}

\makeatletter
\def\author@andify{%
  \nxandlist {\unskip ,\penalty-1 \space\ignorespaces}%
    {\unskip {} \@@and~}%
    {\unskip \penalty-2 \space \@@and~}%
}

\title{PAC structures in nutshell}

\author[D. M. HOFFMANN]{Daniel Max Hoffmann$^{\dagger}$}
\thanks{$^{\dagger}$SDG. The author is supported by the Narodowe Centrum Nauki grants no. 2016/21/N/ST1/01465,
and 2015/19/B/ST1/01150.}
\address{$^{\dagger}$Instytut Matematyki\\
Uniwersytet Warszawski\\
Warszawa\\
Poland}
\email{daniel.max.hoffmann@gmail.com}

\DeclareMathOperator{\fratt}{Fratt}

\DeclareMathOperator{\acl}{acl} \DeclareMathOperator{\dcl}{dcl} 
 \DeclareMathOperator{\aut}{Aut} \DeclareMathOperator{\id}{id}

\DeclareMathOperator{\stab}{Stab}

 \DeclareMathOperator{\theo}{Th}\DeclareMathOperator{\alg}{alg}

 \DeclareMathOperator{\eq}{eq}

\DeclareMathOperator{\ddf}{DF}\DeclareMathOperator{\dcf}{DCF}\DeclareMathOperator{\scf}{SCF}

\DeclareMathOperator{\tcf}{TCF}

\DeclareMathOperator{\mc}{mc}

\DeclareMathOperator{\Cn}{Cn}
\DeclareMathOperator{\qftp}{qftp}

\DeclareMathOperator{\res}{res}
\DeclareMathOperator{\ST}{ST}

\DeclareMathOperator{\ec}{ec}

\newtheorem{theorem}{Theorem}[section]
\newtheorem{prop}[theorem]{Proposition}
\newtheorem{lemma}[theorem]{Lemma}
\newtheorem{cor}[theorem]{Corollary}
\newtheorem{fact}[theorem]{Fact}
\theoremstyle{definition}
\newtheorem{definition}[theorem]{Definition}
\newtheorem{example}[theorem]{Example}

\theoremstyle{remark}
\newtheorem*{theorem*}{Theorem}
\newtheorem*{cor*}{Corollary}

\theoremstyle{definition}
\theoremstyle{definition}

\theoremstyle{definition}

\theoremstyle{remark}

\makeatletter
\providecommand*{\cupdot}{%
  \mathbin{%
    \mathpalette\@cupdot{}%
  }%
}
\newcommand*{\@cupdot}[2]{%
  \ooalign{%
    $\m@th#1\cup$\cr
    \sbox0{$#1\cup$}%
    \dimen@=\ht0 %
    \sbox0{$\m@th#1\cdot$}%
    \advance\dimen@ by -\ht0 %
    \dimen@=.5\dimen@
    \hidewidth\raise\dimen@\box0\hidewidth
  }%
}

\providecommand*{\bigcupdot}{%
  \mathop{%
    \vphantom{\bigcup}%
    \mathpalette\@bigcupdot{}%
  }%
}
\newcommand*{\@bigcupdot}[2]{%
  \ooalign{%
    $\m@th#1\bigcup$\cr
    \sbox0{$#1\bigcup$}%
    \dimen@=\ht0 %
    \advance\dimen@ by -\dp0 %
    \sbox0{\scalebox{2}{$\m@th#1\cdot$}}%
    \advance\dimen@ by -\ht0 %
    \dimen@=.5\dimen@
    \hidewidth\raise\dimen@\box0\hidewidth
  }%
}
\makeatother

\def\Ind#1#2{#1\setbox0=\hbox{$#1x$}\kern\wd0\hbox to 0pt{\hss$#1\mid$\hss}
\lower.9\ht0\hbox to 0pt{\hss$#1\smile$\hss}\kern\wd0}

\def\ind{\mathop{\mathpalette\Ind{}}}

\def\notind#1#2{#1\setbox0=\hbox{$#1x$}\kern\wd0
\hbox to 0pt{\mathchardef\nn=12854\hss$#1\nn$\kern1.4\wd0\hss}
\hbox to 0pt{\hss$#1\mid$\hss}\lower.9\ht0 \hbox to 0pt{\hss$#1\smile$\hss}\kern\wd0}

\begin{document}
\newcommand{\ov}{\overline}
\newcommand{\FC}{\mathfrak{C}}

\newcommand{\twoc}[3]{ {#1} \choose {{#2}|{#3}}}
\newcommand{\thrc}[4]{ {#1} \choose {{#2}|{#3}|{#4}}}
\newcommand{\Rr}{{\mathds{R}}}
\newcommand{\Kk}{{\mathds{K}}}

\newcommand{\dlog}{\mathrm{ld}}
\newcommand{\ga}{\mathbb{G}_{\rm{a}}}
\newcommand{\gm}{\mathbb{G}_{\rm{m}}}
\newcommand{\gaf}{\widehat{\mathbb{G}}_{\rm{a}}}
\newcommand{\gmf}{\widehat{\mathbb{G}}_{\rm{m}}}
\newcommand{\gdf}{\mathfrak{g}-\ddf}
\newcommand{\gdcf}{\mathfrak{g}-\dcf}
\newcommand{\fdf}{F-\ddf}
\newcommand{\fdcf}{F-\dcf}
\newcommand{\mw}{\scf_{\text{MW},e}}

\newcommand{\BC}{{\mathbb C}}

\newcommand{\CC}{{\mathcal C}}
\newcommand{\CG}{{\mathcal G}}
\newcommand{\CK}{{\mathcal K}}
\newcommand{\CL}{{\mathcal L}}
\newcommand{\CN}{{\mathcal N}}
\newcommand{\CS}{{\mathcal S}}
\newcommand{\CU}{{\mathcal U}}
\newcommand{\CF}{{\mathcal F}}
\newcommand{\CP}{{\mathcal P}}
\newcommand{\CI}{{\mathcal I}}

\maketitle
\vspace{-8mm}
\begin{center}
{\small Instytut Matematyki, Uniwersytet Warszawski}
\end{center}

\section{Introduction}
This is an overwiev of recent results about \emph{pseudo-algebraically closed} substructures (PAC substructures) of some stable ambient structure. Our aim is a compact presentation of the main results from \cite{Hoff3}, \cite{Hoff4}, \cite{DHL1} and \cite{HL1}, so a reader interested in studying PAC substructures might quickly get into the subject. 
The content of the text covers all topics from our talk at the RIMS Model Theory Workshop 2018 and the arrangement of the text is based on the order of the mentioned talk.

We fix some stable theory $T$ in a language $\mathcal{L}$, such that $T$ has quantifier elimination and elimination of imaginaries, then we choose a monster model of $T$, say $\mathfrak{C}$, and a group $G$.

We are interested in providing new examples of structures which may serve in the studies on neo-stability. So the general rule for our approach is to start with a stable structure and then ``break" the stability in a controlled manner. It turns out that PAC substructures suit 	perfectly for our purposes, since they are controlled by absolute Galois groups. In other words, we want to ``break" stability and obtain PAC substructures with understandable absolute Galois groups.

We thank organizers of the RIMS Model Theory Workshop 2018 - especially, we are grateful to Hirotaka Kikyo for his support and kindness.

\section{Pseudo-algebraically closed structures}
\subsection{Definitions}
In the theory of fields an extension of fields $E\subseteq F$ is \emph{regular} if it is separable and $E$ is relatively algebraically closed in $F$ ($E^{\alg}\cap F=E$). There is no reasonable generalization of the notion of separability for abstract structures, therefore, if we want to introduce a general notion of regularity, we need to restrict our considerations to some subclass of substructures of $\FC$. If the field $E$ is perfect (i.e. definably closed in some monster model of ACF), then $E\subseteq F$ is separable. Passing to definable closure seems to be a solution:

\begin{definition}\label{regular.def}
Let $E\subseteq A$ be small subsets of $\mathfrak{C}$. We say that $E\subseteq A$ is \emph{regular} if
$$\dcl(A)\cap\acl(E)=\dcl(E).$$
\end{definition}

The above notion of regularity coincides with the notion of stationarity given in \cite[Definition 5.17]{silvain002} and in \cite{hrushovski_onfinima}:
\begin{fact}[Lemma 3.35 in \cite{Hoff3}]\label{PACclaim}
For a small set $E\subseteq\mathfrak{C}$ and a complete type $p$ over $E$ it follows:
\begin{IEEEeqnarray*}{rCl}
p\text{ is stationary} &\iff & (\forall A\models p)(E\subseteq EA\text{ is regular}) \\
 &\iff & (\exists A\models p)(E\subseteq EA\text{ is regular}).
\end{IEEEeqnarray*}
\end{fact}
Note that regularity is invariant under the action of automorphisms, every extension of an algebraically closed subset is regular, and regularity is transitive (if $E\subseteq A$ and $A\subseteq B$ are regular, then also $E\subseteq B$ is regular). By Lemma 3.3 in \cite{Hoff3}, if $A$ is an existentially closed substructure of $B$ (where $A$ and $B$ are small substructures of $\FC$), then $A\subseteq B$ is regular, so regular extensions are related to existential closedness.
A substructure is existentially closed if it is existentially closed in all extensions (i.e. all small extensions living inside $\FC$). Quite interesting notion occurs, if we relax conditions in this definition by requiring existential closedness only in regular extensions:

\begin{definition}\label{def:PAC}
Assume that $M\preceq\FC$.
Let $N$ be a substructure of $M$. We say that $N$ is \emph{pseudo-algebraically closed} (\emph{PAC})
in $M$ if for every small (in the sense of $\FC$) substructure $N'$ of $M$, which is a regular extension of $N$, it follows that $N$ is existentially closed in $N'$. If $M=\FC$ we say ``PAC substructure" instead of ``PAC substructure in $\FC$".
\end{definition}

Note that PAC substructures are definably closed (so PAC substructures in the case of fields coincide with perfect PAC fields).
For some reasons it will be better to call PAC substructures ``regularly closed substructures" (e.g. it may happen that algebraically closed substructure is not pseudo-algebraically closed), but we introduced this definition as an analogon of the well known (in field theory) notion of PAC fields. A field $K$ is pseudo-algebraically closed (is a PAC field) if each non-empty absolutely irreducible $K$-variety has a $K$-rational point. It might be difficult to generalize such a definition of PAC fields, but luckily for us there is a theorem saying that PAC fields are exactly these fields which are existentially closed in every regular extension (Proposition 11.3.5 in \cite{FrJa}). Hence ``regularly closed substructures" are called ``PAC substructures". Introducing PAC substructures (called also PAC structures) was motivated by key role of PAC fields in the research of field theory in the second half of the 20th century. In short,
they were ``discovered" in \cite{Ax1} and \cite{Ax2} as a generalization of pseudofinite fields, but the name ``PAC fields" was given in \cite{FREY}.
Reader interested in results about PAC fields may consult:
\cite{chatzidakis2002}, \cite{chahru04}), \cite{ChaPil}, \cite{cherlindriesmacintyre}, \cite{CDM81}, 
\cite{ershov1980}.
The most important property of PAC fields - so-called ``Elementarily Equivalence Theorem" (Theorem 20.3.3 in \cite{FrJa}, Proposition 33 in \cite{cherlindriesmacintyre}) - implies that, to a large extent, logical properties of PAC fields are controlled by their absolute Galois groups.
A desired generalization - ``Elementarily Equivalence Theorem for Structures" - will be provided in an upcoming subsection (two PAC substructures have the same first order theory if they have isomorphic Galois groups), but before this let us discuss basic properties of PAC structures.

In our approach, we define PAC substructures as a generalization of existentially closed substructures. However, we are not first in studying PAC substructures. Studying PAC substructures has begun (to our knowledge) in \cite{manuscript}, where the author attempted to develop ``geometric simplicity theory" by developing an analogon of model theory of bounded PAC fields (which are simple) in the context of strongly minimal ambient model (in our case the ambient model is stable). He defined PAC substructures in the following way:

\begin{definition}[Definition 1.2 in \cite{manuscript}]\label{def:hru.PAC}
Let $T'$ be strongly minimal and has quantifier elimination and elimination of imaginaries. Let $M\models T'$ satisfy the \emph{definable multiplicity property} (consult ``Framework" at page 8. in \cite{manuscript}). A subset $P$ of $M$ is a \emph{PAC$_{\text{Hru}}$} subset of $M$, if every multiplicity $1$ formula with parameters from $P$ has a solution in $P$.
\end{definition}

By Proposition 3.10 in \cite{Hoff3}, both notions of a PAC substructure (Definitions \ref{def:PAC} and \ref{def:hru.PAC}) coincide in the strongly minimal context. Many important facts about PAC structures are provided in \cite{PilPol} and in \cite{Polkowska}. Authors of \cite{PilPol} and \cite{Polkowska} work with a stronger notion of a PAC substructure (in the stable context):

\begin{definition}[Definition 3.1 in \cite{PilPol}]\label{PilPol.PAC.def}
Let $T'$ be stable, $\kappa\geqslant |T'|^+$ a cardinal, $M\models T'$, 
and let $P$ be an $\mathcal{L}$-substructure of $M$. We say that $P$ is a \emph{$\kappa$-PAC$_{\text{PP}}$} 
substructure of $M$ if whenever $A\subseteq P$ has cardinality smaller than $\kappa$ and $p(x)$ is complete stationary type over $A$ (in the sense of $M$), then $p$ has a realization in $P$.
\end{definition}
\noindent
This definition of a PAC substructure assumes some saturation, which - as we will see soon - is crucial if we want to avoid serious obstructions (e.g. in the Elementary Equivalence Theorem for Structures). For the reader convenience we visualize relations between all three definitions.
In the set-up of a strongly minimal theory:
$$\text{PAC}_{\text{Hru}}\;\;+\;\;\kappa-\text{saturated}\quad\Rightarrow\quad\kappa-\text{PAC}_{\text{PP}}\quad\Rightarrow\quad\text{PAC}_{\text{Hru}}=\text{PAC}.$$
In the set-up of a stable theory:
$$\text{PAC}\;\;+\;\;\kappa-\text{saturated}\quad\Rightarrow\quad\kappa-\text{PAC}_{\text{PP}}\quad\Rightarrow\quad\text{PAC}.$$
It was mentioned that working with somehow saturated PAC structures might be more convenient. Nevertheless, passing to an elementary extension which is PAC substructure is not ``free of charge": 

\begin{lemma}[Lemma 2.12 in \cite{DHL1}]\label{lemma:saturated_PAC_str}
Suppose PAC is a first order property and 
pure saturation over $P$ is a first order property. 
Let $P$ be PAC in $\mathfrak{C}$. Let $\kappa$ be an infinite cardinal. Then, there is a $\kappa$-saturated elementary extension $P^*$ of $P$, which is PAC in $\mathfrak{C}$. 
\end{lemma}
\noindent
Where the missing notions are defined as follows:

\begin{definition}\label{def:PAC2}
Assume that $M\preceq\mathfrak{C}$ and $P$ is a substructure of $M$. 
\begin{enumerate}
\item We say that $M$ is \emph{purely saturated over $P$} if every type over $P$ is realized in $M$.

\item We say that \emph{PAC is a first order property} if there exists a set $\Sigma$ of $\mathcal{L}\cup\{P\}$-formulas such that
$$(M,P)\models T\cup\Sigma\qquad\iff\qquad M\models T\text{ and }P\text{ is PAC in }M.$$

\item We say that \emph{pure saturation over $P$ is a first order property} if there exists a set $\Sigma$ of $\mathcal{L}\cup\{P\}$-formulas such that
\begin{enumerate}
\item[(i)] if $M\models T$ is purely saturated over $P$, then $(M,P)\models T\cup\Sigma$,
\item[(ii)] if $(M,P)\models T\cup\Sigma$ and $(M,P)$ is $|T|^+$-saturated, then $M$ is purely saturated over $P$.
\end{enumerate}
\end{enumerate}
\end{definition}
\noindent
If $T$ has no finite cover property
(an assumption stronger than stability, see Definition 4.1 and Theorem 4.2 in Chapter II of \cite{ShelahModels}), then pure saturation over $P$ is a first order property (see Remark 2.10 in \cite{DHL1} or Remark 3.6 in \cite{PilPol}). Therefore our impression is that it might be very hard to pass to a saturated elementary extension which remains PAC (which is not the case for existentially closed substructures) and so results for non-saturated PAC substructures are desired, even if they require more complicated techniques.
Some of results require a weaker variant of saturation which we introduce right now.

For a small substructure $P$ of $\mathfrak{C}$, let us define 
\begin{IEEEeqnarray*}{rCl}
\ST(P,\kappa,\lambda) &:=& \{\qftp(\bar{d}/A)\;\;|\;\;A\subseteq P,\;|A|<\kappa, \\
 & & P\subseteq\dcl(P,\bar{d})\text{ is regular and }|\bar{d}|<\lambda\}.
\end{IEEEeqnarray*}

If $P$ is PAC, then each element of $\ST(F,\kappa,\lambda)$ is a partial type in the sense of $\theo(P)$. We say that $P$ is 
 \emph{$_{\kappa}^{\lambda}$-regularly} saturated if every element of $\ST(P,\kappa,\lambda)$ is realized in $P$. We use ``$\kappa$-regularly saturated" for ``$^{\kappa}_{\kappa}$-regularly saturated" (see Definition 3.2 in \cite{DHL1}). 
If $P$ is PAC and $\kappa$-saturated (in the sense of $\theo(P)$), then $P$ is $^{\omega}_{\kappa}$-saturated.

On the other hand passing to elementary substructure preserves being PAC (by Fact 2.4 and Fact 2.6 in \cite{DHL1}).

\subsection{A bite of Galois theory}
First of all, let us clarify what do we understand by \emph{Galois extension}.

\begin{definition}\label{galois.ext.def}
\begin{enumerate}
\item Assume that $A\subseteq C$ are substructures of $\mathfrak{C}$. We say that $C$ is \emph{normal over $A$} (or we say that $A\subseteq C$ is a \emph{normal extension}) if $\aut(\mathfrak{C}/A)\cdot C\subseteq C$.
(Note that if $C$ is small and $A\subseteq C$ is normal, then it must be $C\subseteq\acl(A)$.)

\item Assume that $A\subseteq C\subseteq\acl(A)$ are small substructures of $\mathfrak{C}$ such that $A=\dcl(A)$, $C=\dcl(C)$ and $C$ is normal over $A$. In this situation we say that $A\subseteq C$ is a \emph{Galois extension}.
\end{enumerate}
\end{definition}

If $A\subseteq C$ is a Galois extension, then $\aut(C/A)$ is a profinite group.
Many facts from classical Galois theory have their counterparts in the model-theoretic variant of Galois theory. Some of these generalization might be found
in \cite{invitation},
in Section 3. in \cite{Hoff3} or in \cite{Hoff4}. To convince the reader that a lot of the classical Galois theory generalizes quite straightforward, let us provide some highlights.

\begin{fact}[The Galois correspondence]\label{galois.correspondence}
Let $A\subseteq C$ be a Galois extension, introduce
$$\mathcal{B}:=\lbrace B\;|\;A\subseteq B=\dcl(B)\subseteq C\rbrace,$$
$$\mathcal{H}:=\lbrace H\;|\; H\leqslant\aut(C/A)\text{ is closed}\rbrace.$$
Then $\alpha(B):=\aut(C/B)$ is a mapping between $\mathcal{B}$ and $\mathcal{H}$, $\beta(H):=C^H$ is a mapping between $\mathcal{H}$ and $\mathcal{B}$ and it follows
$$\alpha\circ\beta=\id,\qquad\beta\circ\alpha=\id.$$
\end{fact}

\begin{lemma}[Lemma 2.10 in \cite{Hoff4}]\label{N_Galois}
Assume that $N$ is a small definably closed substructure of $\mathfrak{C}$ equipped with a $G$-action $(\tau_g)_{g\in G}$ ($|G|$ is smaller than the saturation of $\mathfrak{C}$). 
Let $i:G\to\aut(N/N^G)$ be given by $i(g):=\tau_g$.
\begin{enumerate}
\item If for every $b\in N$ the orbit $G\cdot b$ is definable, then $N^G\subseteq N$ is normal.

\item
If $N\subseteq\acl(N^G)$, then $N^G\subseteq N$ is a Galois extension.

\item
If $G$ is finite, then $N\subseteq\acl(N^G)$, hence also the second point follows.

\item (Artin's theorem)
If $G$ is finite and the $G$-action $(\tau_g)_{g\in G}$ is faithful, then 
$i:G \cong \aut(N/N^G)$.

\item
Assume that $G$ is profinite, the $G$-action $(\tau_g)_{g\in G}$ is faithful and for every $m\in N$ the stabiliser $\stab(m)=\lbrace g\in G\;|\;\tau_g(m)=m\rbrace$ is an open subgroup of $G$. Then $N^G\subseteq N$ is Galois and
$i: G\cong \aut(N/N^G)$ (as profinite groups).
\end{enumerate}
\end{lemma}

\emph{Absolute Galois groups} and \emph{Galois groups of extensions} play a major role in our description of PAC substructures of $\mathfrak{C}$:
\begin{IEEEeqnarray*}{rCl}
G(A):=\aut(\acl(A)/\dcl(A)) &\qquad& \text{(the absolute Galois group of $A$)}\\
G(B/A):=\aut(\dcl(B)/\dcl(A)) &\qquad &\text{(the Galois group of extension $A\subseteq B$)}
\end{IEEEeqnarray*}
For example, $A\subseteq B$ is regular extension if and only if the restriction map $G(B)\to G(A)$ is onto.
Moreover, one may find interesting that all profinite groups occur in our stable theory $T$ as Galois groups of some Galois extensions:

\begin{fact}[Corollary 3.3 in \cite{Hoff4}]\label{pro.embed2}
A group $G$ is profinite if and only if there exist a Galois extension $A\subseteq B$ of small substructures of $\mathfrak{C}$ such that $G\cong\aut(B/A)$.
\end{fact}

An important property of the absolute Galois group of a PAC substructure is projectivity, which was known in the strongly minimal context (Lemma 1.17 in \cite{manuscript}) and then generalized to the stable context (Theorem 4.4 in \cite{Hoff4}):

\begin{theorem}[Theorem 4.4 in \cite{Hoff4}]\label{PAC.proj}
If a small $N$ is PAC, then $G(N)$ is projective.
\end{theorem}

A natural question is whether any projective profinite group occurs as the absolute Galois group of some PAC substructure. The following result puts some light on this problem.

\begin{theorem}\label{fr.ja.2311}
Assume that $A\subseteq B$ is a Galois extension of small substructures of $\mathfrak{C}$ and assume that there is an epimorhism of profinite groups
$\alpha:G\to\aut(B/A)$, and $G$ is projective.
there exists a definably closed substructure $P\supseteq A$ of $\mathfrak{C}$ 
and an isomorphism of profinite groups $\gamma:G\to G(P)$ such that
$$\xymatrix{G \ar[r]^-{\gamma} \ar[dr]_-{\alpha} &  G(P) \ar[d]^{|_N} \\
& \aut(B/A)
}$$
is commuting.
Moreover, if any type over $A$ has only finitely many extensions over $\acl(A)$, then $P$ is PAC.
\end{theorem}

For example, if $T$ is $\omega$-stable, then a profinite group $G$ is projective if and only if $G$ is isomorphic to the absolute Galois group of some PAC substructure of $\mathfrak{C}$ (see Corollary 4.10 in \cite{Hoff4}). How about the same for any stable $T$? We do not know it due to the use of a generalization of the Ax-Roquette theorem in the proof of Theorem \ref{fr.ja.2311}.
There are problems with generalizing the Ax-Roquette theorem (algebraic extension of a PAC field is a PAC field), more precisely: there is an example of a stable theory $T$ where algebraic closure of a PAC substructure is not a PAC substructure (see Section 5. in \cite{PilPol}). The best (to our knowledge) generalization of the Ax-Roquette theorem was provided in Proposition 3.9 in \cite{PilPol}:

\begin{prop}
Suppose $P$ is $\kappa$-PAC$_{\text{PP}}$ substructure of $\mathfrak{C}$.
\begin{enumerate}
\item[i)] Any finite algebraic extension of $P$ is also a $\kappa$-PAC$_{\text{PP}}$ substructure of $\mathfrak{C}$.
\item[ii)] Suppose that in $T$ every (finitary) complete type over any set has finite
multiplicity. Then any algebraic extension of $P$ is $\kappa$-PAC$_{\text{PP}}$.
\end{enumerate}
\end{prop}

\noindent
In Lemma 4.5 in \cite{Hoff4}, the above proposition was formulated for PAC substructures: if $P$ is a PAC substructure and any type over $P$ has only finitely many (non-forking) extensions over $\acl(P)$, then any algebraic extension of $P$ is a PAC substructure. Actually, formulation of Lemma 4.5 in \cite{Hoff4}, which is rather complicated, gives us also that any finite Galois extension of a PAC substructure (``finite" means that the Galois group of this extension is finite) is a PAC substructure (to show this you may use ``Primitive Element Theorem", Theorem 5.3 in \cite{DHL1}).

\subsection{Elementarily Equivalence Theorem}
The idea that the theory of a PAC substructure might be controlled by its absolute Galois group comes from the theory of PAC fields (see works of Fried, Jarden and Kiehne). 
However, some results for PAC substructures were obtained even without generalizing the Embedding Lemma (Lemma 2.1 in \cite{JardenKiehne}) or Elementarily Equivalence Theorem (Theorem 3.2 in \cite{JardenKiehne}). For example, in \cite{Polkowska}, it was shown that PAC substructures with small absolute Galois group (i.e. bounded PAC substructures) are simple (Corollary 3.22 in \cite{Polkowska}). This result was obtained under assumption that ``PAC is a first order property", which seems to be natural assumption in studying PAC substructures.

The Embedding Lemma was crucial in the proof of the Elementarily Equivalence Theorem for PAC fields and its generalized form (Lemma 3.5 in \cite{DHL1}) plays important role in the proof of the Elementarily Equivalence Theorem for PAC substructures (Proposition 3.8 in \cite{DHL1}). Moreover, the generalized Embedding Lemma might be used in a description of types in PAC substructures and is used in the proof of the upcoming Weak Indpendence Theorem (\cite{HL1}). Since the generalized Embedding Lemma is involved in several techniques, we provide it here in the most general formulation.

\begin{lemma}[Embedding Lemma, Lemma 3.5 in \cite{DHL1}]\label{lemma2022}
Assume that
\begin{itemize}
\item $L\subseteq L'$, $M\subseteq M'$, $E\subseteq E'$, $F\subseteq F'$ are small Galois extensions in $\mathfrak{C}$,

\item $L\subseteq E$, $M\subseteq F$, $M'\subseteq F'$,
\item $L'\subseteq E'$ is regular,

\item $F$ is $\kappa$-regularly saturated, where $\kappa\geqslant(\max\{|E|,|T|\})^+$,
\item $F$ is PAC in $\mathfrak{C}$,
\item $\Phi_0\in\aut(\mathfrak{C})$ is such that $\Phi_0(L)=M$ and $\Phi_0(L')=M'$,
\item $\varphi: G(F'/F)\to G(E'/E)$ is a continuous group homomorphism such that
$$\xymatrix{
 G(F'/F) \ar[r]^-{\varphi} \ar[d]_{\res} &  G(E'/E) \ar[d]^{\res} \\
 G(M'/M) \ar[r]_{\varphi_0} &  G(L'/L)
}$$
where $\varphi_0(\sigma):=\Phi_0^{-1}\circ \sigma\circ\Phi_0$ for each $\sigma\in G(F'/F)$, commutes.
\end{itemize}
Then there exists $\Phi\in\aut(\mathfrak{C})$ such that
\begin{itemize}
\item $\Phi|_{L'}=\Phi_0|_{L'}$,
\item $\Phi(E)\subseteq F$, $\Phi(E')\subseteq F'$,
\item $\varphi(\sigma)=\Phi^{-1}\circ \sigma\circ \Phi$ for any $\sigma\in G(F'/F)$ .
\end{itemize}
Moreover, if $\varphi$ is onto and $E'=\acl(E)$, then $\Phi(E)\subseteq F$ is regular.
\end{lemma}

On the other hand, we provide here only a consequence, which has rather simple and elegant formulation, from the Elementarily Equivalence Theorem for PAC substructures (Proposition 3.8 in \cite{DHL1}).

\begin{theorem}[Elementarily Equivalence Theorem, Corollary 3.9 in \cite{DHL1}]\label{cor2033}
If $E$ and $F$ are $\kappa$-regularly saturated PAC substructures of $\mathfrak{C}$ ($\kappa\geqslant|T|$),
and for some definably closed $L\subseteq F\cap E$ of size strictly smaller than $\kappa$ there exists a continuous isomorphism $\varphi: G(F)\to G(E)$ such that
$$\xymatrix{ G(F) \ar[dr]_{\res} \ar[rr]^{\varphi}& &  G(E)\ar[dl]^{\res} \\
&  G(L) &
}$$
is commuting, then $E\equiv_L F$.
\end{theorem}

It is remarkable that there is no (to our knowledge) easy way to relax the assumption about saturation in the above Theorem \ref{cor2033}. We can assume that ``PAC is a first order property" and ``pure saturation over $P$ is a first order property" and then try to pass to saturated elementary extensions of our given PAC substructures $E$ and $F$. The problem in this approach is that the isomorphism $\varphi$ of profinite groups $G(E)$ and $G(F)$ does not carry enough data about ``model-theoretic" structure of $E$ and $F$, and therefore there are serious problems in lifting up $\varphi$ to the level of saturated elementary extensions of $E$ and $F$. The issue shows up if we consider many sorted structures:

\begin{example}
Consider two sorted structure $(\mathbb{C},\mathbb{C})$, where $\mathbb{C}$ is the
field of complex numbers $\mathbb{C}$ and there is no interaction between sorts. Note that $(\mathbb{Q},\acl(\mathbb{Q}))\not\equiv(\acl(\mathbb{Q}),\mathbb{Q})$, but
$$ G(\mathbb{Q},\acl(\mathbb{Q}))\cong G(\mathbb{Q})\times G(\acl(\mathbb{Q})) \cong G(\acl(\mathbb{Q}))\times G(\mathbb{Q}) \cong G(\acl(\mathbb{Q}),\mathbb{Q}).$$
We may even pass to $(\mathbb{C},\mathbb{C})^{\eq}$ to get elimination of imaginaries, and then do the Morleyisation, to obtain quantifier elimination, but absolute Galois groups will remain isomorphic despite substructures do not share the same first order theory.
\end{example}

Clearly, the above example is based on fact that the isomorphism between absolute Galois groups is ``mixing" sorts. To avoid such phenomena, the notion of \emph{sorted isomorphism} is introduced (in \cite{DHL1}) and developed into notion of \emph{sorted profinite groups} (in \cite{HL1}), which form a category. In short: sorted structure on $G(F)$ codes on which sorts live primitive elements of finite Galois extensions of $F$. Sorted profinite groups can be presented as first order structures in some abstract many sorted language, so called \emph{(sorted) compelete systems}, and, what is more surprising, interpreted in $(\FC,F)$ (structure $\FC$ considered with a predicate for set $F$). The reader interested in more details may find them in \cite{HL1} (which will be published on arxive soon).

The Elemetarily Equivalence Theorem for non-saturated PAC structures, Theorem 5.8 in \cite{DHL1}, assumes ``PAC is a first order property" and ``pure saturation over $P$ is a first order property", and requires that isomorphism between $G(F)$ and $G(E)$ preserves sorts. So the price for relaxing the assumption about saturation (in Thereom \ref{cor2033}) is high and to avoid paying it, in the following text, we will focus mostly on saturated PAC substructures.

\subsection{Weak Independence Theorem}
This subsection is only to mention some ideas and describe results without precise formulation of them. 
Elementarily Equivalence Theorem characterizes the first order theory of a PAC substructure/field.
Usually model theory aspires to make one more step and try to find some geometrical behavior in models of a given theory. We would like to imitate methods from geometric stability theory in non-stable theories which are lying not to far from stability. To do this we need to consider a notion of independence, which inherits some properties from well-known geometrical notions of independence.

In \cite{chatzidakis2017}, the author describes a very interesting concept of inducing a notion of independence in a (saturated) PAC field $K$ by combining forking independence present in the separable closure of $K$ and forking independence given in the complete system corresponding to $G(K)$. The main ingredient is Theorem 2.1 from \cite{chatzidakis2017}, which says that if there is ``a solution of the independence theorem problem" on the level of complete system corresponding to $G(K)$, then there is ``a solution of the independence theorem problem" in the field $K$. Using this result, Nick Ramsey in his doctoral dissertation (Corollary 7.2.7 in \cite{NickThesis}) shows that for a PAC field $K$ the following holds: if the theory of the complete system related to $G(K)$ is NSOP$_1$, then $K$ is NSOP$_1$. In Ramsey's result the notion of independence in $K$ is combined from forking independence in the separable closure of $K$ and from Kim-dividing on the level of complete system corresponding to $G(K)$.

We were interested in generalizing the aforementioned results 
to the level of PAC substructures of stable monster $\FC$
and we plan to provide such generalizations in \cite{HL1}. For example, we already generalized Theorem 2.1 from \cite{chatzidakis2017}, which is the main ingredient also on the level of PAC substructures. Now, generalizing Theorem 7.2.6 and Corollary 7.2.7 from \cite{NickThesis} should be straightforward and we are working on it.

One of our goals is to provide new examples of theories with NSOP$_1$. One could ask: could not we construct such examples in a more direct way, e.g. structures given by more combinatorial approach?
Yes, we could, but in the case of NSOP$_1$ PAC substructures, we already dispose a tool to study model-theoretical properties of these structures - the general Galois theory. And as we will see in a moment, there is a way to construct PAC substructures with a desired absolute Galois group.

\section{Substructures with $G$-action}
This part borrows results from \cite{Hoff3} and \cite{Hoff4}. The reader should treat \cite{Hoff3} as a published version of author's doctoral dissertation. The initial idea behind \cite{Hoff3} was to obtain a way of ``breaking", in a controlled manner, stability of a given stable theory $T$ by considering additional group action by automorphisms. The desired goals were achieved and under some assumptions the resulting theory remains simple but not stable. As side product of studies on existentially closed substructures with a group action, it was observed that such substructures are PAC substructures and the definition of PAC substructures was slightly reformulated (see Section 3.1 in \cite{Hoff3}) so we obtained Definition \ref{def:PAC}. 
It was also observed that there is a link between group $G$ which acts on considered substructures
and their absolute Galois groups. Therefore we may adjust $G$ to control PAC substructures given as reducts of existentially closed substructures with a group action of $G$.

\subsection{Basics}
The research in \cite{Hoff3}, which is a base for the following lines, was motivated by results of \cite{ChaPil} and \cite{nacfa}.
We fix a group $G$ and set language $\mathcal{L}^{G}$ to be the language $\mathcal{L}$ extended by unary function symbols $\bar{\sigma}=(\sigma_g)_{g\in G}$. 

\begin{definition}[Definition 2.5 in \cite{Hoff3}]
\begin{enumerate}
\item We introduce set of $\mathcal{L}^G$-sentences $A_G$, which contains exactly the following axioms:
\begin{enumerate}
\item[i)] $\sigma_g$ is an automorphism of $\mathcal{L}$-structure for every $g\in G$,
\item[ii)] $\sigma_g\circ\sigma_h=\sigma_{g\ast h}$ for every $g,h\in G$.
\end{enumerate}
\item
Let $(M,(\sigma_g)_{g\in G})$ be an $\mathcal{L}^G$-structure. We say that $(\sigma_g)_{g\in G}$ is a \emph{$G$-action} on $M$ if $(M,(\sigma_g)_{g\in G})\models A_G$.

\item
If $T'$ is an $\mathcal{L}$-theory, then $T'_G$ is an $\mathcal{L}^G$-theory
equal to the set of consequences of $T'\cup A_G$, i.e. $T'_G=\Cn(T'\cup A_G)$.
\end{enumerate}
\end{definition}
There is an interesting question whether a class of existentially closed models of $T'_G$, for a given theory $T'$, is an elementary class, in other words whether model companion of $T'_G$ exists (for some reasons it is more reasonable to ask about the existence of the model companion for $(T'_{\forall})_G$).
The question is even more interesting, since the answer for the case of $G=\mathbb{Z}$ is related to nfcp.
For research in this line, the reader may consult \cite{balshe}, \cite{kikyo1}, \cite{kikyo2}, \cite{kikpil} and \cite{kishe}.

For our purposes, we will focus on the class of existentially closed models of $(T_{\forall})_G$ (recall that $T$ is stable and allows to eliminate quantifiers and imaginaries), so small substructures of $\FC$ which are equipped with a $G$-action and which are existentially closed among the class of substructures equipped with a $G$-action.
By Proposition 3.56 in \cite{Hoff4}, such substructures are PAC substructures in $\FC$. Moreover, if $G$ is finitely generated, then also the substructure of invariants is PAC (by Proposition 3.51 in \cite{Hoff3}) and it is bounded (Corollary 3.50 in \cite{Hoff3}). More precisely, if $(M,\bar{\sigma})$ is an existentially closed model of $(T_{\forall})_G$, then by invariants of the $G$-action we understand elements of the set $M^G:=\{m\in M\;|\; \sigma_g(m)=m\text{ for all }g\in G\}$. By Corollary 3.22 in \cite{Polkowska}, the set $M^G$ for finitely generated $G$, constitutes an $\mathcal{L}$-substructure which is, as bounded PAC substructure, simple (see Theorem 4.40 in \cite{Hoff3}). Now, we proceed to the simplicity of a prospective model companion of the theory $(T_{\forall})_G$.

\subsection{Independence Theorem and simplicity}
In this subsection we assume that the model companion of the theory $(T_{\forall})_G$, say $T_G^{\mc}$, exists. 
We will describe assumptions needed to get simplicity of $T_G^{\mc}$.
First, let us recall that if a PAC substructure is bounded (its absolute Galois group is small as a profinite group) then its theory is simple (Corollary 3.22 in \cite{Polkowska}). In the case of PAC fields, we have a stronger fact: a PAC field $K$ is simple if and only if it is bounded (Fact 2.6.7 in \cite{kim1}).
If $(K,\bar{\sigma})$ is an existentially closed field with a $G$-action, then $K$ is PAC and if $K$ is not bounded then $(K,\bar{\sigma})$ can not be simple (since it reduct is not simple). Therefore, if we attempt to show simplicity of some monster model $(\mathfrak{M},\bar{\sigma})$ of $T_G^{\mc}$, it is natural to assume that $\mathfrak{M}$ is bounded in $\FC$ (i.e. $G(\mathfrak{M})$ is a small profinite group, e.g. if $G$ is finite then $\mathfrak{M}$ is bounded, see Proposition 4.25 in \cite{Hoff3}):

\begin{theorem}[Theorem 4.21 and Proposition 4.26 in \cite{Hoff3}]
Assume that $\mathfrak{M}$ is bounded in $\FC$.
Let $(M,\bar{\sigma})\preceq(\mathfrak{M},\bar{\sigma})$ and let $p_1(x_1)$, $p_2(x_2)$, $p_3(x_3)$, $p_{12}(x_1,x_2)$, $p_{23}(x_2,x_3)$ and $p_{13}(x_1,x_3)$ be complete $\mathcal{L}^G$-types over $M$ which satisfy $p_i(x_i),p_j(x_j)\subseteq p_{ij}(x_i,x_j)$ and if $a_ia_j\models p_{ij}(x_i,x_j)$ then
$$a_j\ind_M^{\circ} a_i.$$
Then, there exists a complete $\mathcal{L}^G$-type $p_{123}(x_1,x_2,x_3)$ which extends each $p_{ij}(x_i,x_j)$ and such that if $a_1a_2a_3\models p_{123}(x_1,x_2,x_3)$ then
$$a_3\ind_M^{\circ}a_1a_2.$$
\end{theorem}

Ternary relation $\ind^{\circ}$ which occurs in  the above theorem is defined as follows:
$$A\ind^{\circ}_E B\quad\iff\quad G\cdot A\ind_{G\cdot E}G\cdot B,$$
where $A,B,E$ are small subsets of $\mathfrak{M}$. By Proposition 4.16 from \cite{Hoff3}, this ternary relation is an independence relation, hence it is forking independence in $(\mathfrak{M},\bar{\sigma})$ and so $(\mathfrak{M},\bar{\sigma})$ is simple (Theorem 4.22 in \cite{Hoff3}) - under the assumption that $\mathfrak{M}$ is bounded. Moreover, boundedness of $\mathfrak{M}$ implies also that the theory of $(\mathfrak{M},\bar{\sigma})$ has geometric elimination of imaginaries (Theorem 4.36 in \cite{Hoff3}). 

There are several examples of $(T_{\forall})_G$ which have model companions:
\begin{itemize}
\item ACFA is supersimple, not stable,
\item if $G$ is finite then $G-\tcf$ (e.c. fields with group action of $G$, see \cite{nacfa}) is supersimple of SU-rank equal to $|G|$, not stable
\item $\mathbb{Q}$ACFA (e.c. fields with group action of $\mathbb{Q}$, see \cite{qacfa}) is simple, not supersimple, not stable
\item CCMA (Complex Compact Manifolds with Automorphism, see \cite{hils0}) is super\-simple, not stable
\end{itemize}

We end this part with an unexpected result binding together simplicity of a field equipped with a group action with simplicity of its reduct to the structure of the pure field.

\begin{cor}[Corollary 4.31 in \cite{Hoff3}]\label{fields.cor}
Let $T$ be the theory of fields, $\mathcal{L}$ be the language of rings
and let $G$ be a group (not necessarily finite).
Assume that $T_G^{\mc}$ (= $G-\tcf$) exists and let $(K,(\sigma_g)_{g\in G})\models T_G^{\mc}$. The following are equivalent.
\begin{enumerate}
\item
The theory of $K$ in the language $\mathcal{L}$ is simple.

\item
The field $K$ is a bounded field.

\item
The theory of $(K,(\sigma_g)_{g\in G})$ in the language $\mathcal{L}^G$ is simple.
\end{enumerate}
\end{cor}

\subsection{Frattini cover}
As we already mentioned, we can adjust absolute Galois groups (in the sense of $\FC$) of existentially closed substructures of $\FC$ equipped with a $G$-action by changing the group $G$. Assume that $(M,\bar{\sigma})$ is such a substructure of $\FC$. Now, we will describe some correlations between $G$ and Galois groups related to substructure $M$. This description involves the following notion:

\begin{definition}\label{frattini0}
Let $H,H'$ be profinite groups and $\pi: H\to H'$ be a continuous epimorphism.
The mapping $\pi$ is called a \emph{Frattini cover} if for each closed subgroup $H_0$ of $H$, the condition $\pi(H_0)=H'$ implies that $H_0=H$.
If moreover the profinite group $H$ is projective, the map $\pi$ is the unique
\emph{universal Frattini cover} (consult Proposition 22.6.1 in \cite{FrJa}).
\end{definition}

The intuition behind the following results is that the $G$-action on $M$ is determined by its restriction to $F:=\acl(M^G)\cap M$. Note that $M^G\subseteq F$ is a Galois extension (by Lemma 3.28 in \cite{Hoff3}), hence $G(F/M^G)$ is a profinite group and it makes sense to study this group by standard for profinite groups methods. For example, $G(F/M^G)$ is generated as a profinite group by $(\sigma_g|_F)_{g\in G}$ (Proposition 3.30 in \cite{Hoff3}), which implies that if $G$ is finitely generated, then also $G(M^G)$ is finitely generated (as a profinite group, see Proposition 3.49 in \cite{Hoff3}). The following result generalizes a similar one for fields with group action from \cite{sjogren}, the proof in \cite{Hoff4} uses Theorem 1.1 from \cite{NikoSega} which is based on the classification of finite simple groups.

\begin{prop}[Proposition 5.1 in \cite{Hoff4}]\label{profinite.completion}
If $G$ is finitely generated then $G(F/M^G)\cong \hat{G}$ (the profinite completion of $G$).
\end{prop}

Let us assume that $G$ is finitely generated.
We know that $M^G$ and $M$ are PAC substructures of $\FC$. After using some results from \cite{Hoff4} and Theorem \ref{PAC.proj}, we may conclude:

\begin{cor}[Corollary 5.5 in \cite{Hoff4}]
The restriction map 
$$G(M^G)\to G(M\cap\acl(M^G)/M^G)$$
is the universal Frattini cover.
\end{cor}
\noindent
By Proposition \ref{profinite.completion} we obtain the following short exact sequence
$$G\big(M\cap\acl(M^G)\big) \to\mathcal{G}\big(M^G\big)\to G\big(M\cap\acl(M^G)/M^G\big)\cong\hat{G}$$
which gives us a generalization of a result which was known in the case of fields
(\cite{sjogren}, \cite{nacfa}):

\begin{cor}[Corollary 5.6 in \cite{Hoff4}]\label{cor.Galois.descr}
It follows that
\begin{enumerate}
\item
$G(M^G)\cong\fratt(\hat{G})$ (the domain of the universal Frattini cover of the group $\hat{G}$),

\item
$G\big(M\cap\acl(M^G)\big)\cong\ker\Big(\fratt\big(\hat{G}\big)\to \hat{G}\Big)$.
\end{enumerate}
\end{cor}
\noindent
Finally, we see that $M^G$ is a PAC substructure of $\FC$ controlled by its absolute Galois group which is isomorphic to $\fratt(\hat{G})$. If $G$ is finite, then $M\subseteq\acl(M^G)$ (by Lemma \ref{N_Galois}.(2)) so $M\cap\acl(M^G)=M$ and $M$ is a PAC substructure of $\FC$ controlled by its absolute Galois group given as $\ker\Big(\fratt\big(G\big)\to G\Big)$.

\bibliographystyle{plain}
\bibliography{1nacfa2}

\end{document}